\documentstyle[12pt,amscd,amssymb]{amsart}

\newcommand{\QT}{{\Bbb Q}[[t]]}

\newcommand{\hH}{\hat{H}}
\newcommand{\Dws}{D^{(w,\Sigma)}}
\newcommand{\M}{M^{\text{\scriptsize odd}}_{\Sigma}}
\newcommand{\Y}{\Sigma \times {{\Bbb S}}^1}
\newcommand{\X}{X^o}
\newcommand{\D}{D^o}
\newcommand{\B}{B^o}

\newcommand{\inc}{\hookrightarrow}
\newcommand{\ar}{\rightarrow}
\newcommand{\bd}{\partial}
\newcommand{\x}{\times}
\newcommand{\iso}{\cong}
\newcommand{\isom}{\stackrel{\sim}{\ar}}

\newcommand{\pt}{\text{pt}}
\newcommand{\CP}{{\Bbb C \Bbb P}}

\newcommand{\Diff}{\text{Diff}}

\newcommand{\Sym}[1]{\text{Sym}^{#1}}

\newcommand{\PD}{\text{P.D.}}

\newcommand{\cH}{{\cal H}}

\newcommand{\cV}{{\cal V}}

\renewcommand{\AA}{{\Bbb A}}

\newcommand{\DD}{{\Bbb D}}

\newcommand{\QQ}{{\Bbb Q}}

\newcommand{\SS}{{\Bbb S}}
\newcommand{\ZZ}{{\Bbb Z}}

\renewcommand{\a}{\alpha}
\renewcommand{\b}{\beta}

\newcommand{\g}{\gamma}

\newcommand{\f}{\epsilon}

\newcommand{\k}{\kappa}
\newcommand{\s}{\sigma}

\newcommand{\p}{\phi}

\renewcommand{\S}{\Sigma}
\newcommand{\De}{\Delta}

\theoremstyle{plain}
\begingroup
\newtheorem{thm}{Theorem}
\newtheorem{cor}[thm]{Corollary}
\newtheorem{lem}[thm]{Lemma}
\newtheorem{prop}[thm]{Proposition}
\newtheorem{conj}[thm]{Conjecture}
\endgroup

\theoremstyle{definition}
\newtheorem{defn}[thm]{Definition}

\theoremstyle{remark}
\newtheorem{rem}[thm]{Remark}

\title[Donaldson invariants for connected sums]{Gluing formulae for Donaldson 
invariants for connected sums along surfaces}

\author{Vicente Mu\~noz}
\address{Departamento de \'Albegra, Geometr\'{\i}a y Topolog\'{\i}a  \\
Facultad de Ciencias \\
Universidad de M\'alaga, AP.\ 59 \\ 29080 M\'alaga \\ Spain}
\email{vmunoz@@agt.cie.uma.es}
\date{January, 1997}

\begin{document}

\maketitle

\begin{abstract}
  Following our work in~\cite{paper}, we prove a gluing formula for the
  Donaldson invariants of the connected sum of two four-manifolds along 
  surfaces of the same genus $g$,
  self-intersection zero and representing odd homology classes, solving
  a conjecture of Morgan and Szab\'o~\cite{MSz}.
\end{abstract}

\section{Introduction}
  \label{sec:intro}
  This paper tries to answer the question of the behaviour of the
  Donaldson invariants under connected sums along surfaces of 
  arbitrary genus $g \geq 2$ and
  self-intersection zero. The problem was first motivated by the computation
  of the basic classes of elliptic surfaces. The general elliptic surface
  with $b^+>1$ and simply connected, is constructed by the process of 
  connected sum along embedded tori of self-intersection zero from
  elementary pieces (rational elliptic surfaces and homotopy $K3$-surfaces).
  After establishing the appropriate gluing formula for the 
  behaviour of Donaldson invariants under the operation
  of connected sum along embedded tori,
  the information of the basic classes for some particular examples gives 
  the basic classes of any such elliptic surface. This has been carried out
  by many 
  authors~\cite{Morgan}~\cite{Do1}~\cite{Thesis}~\cite{MSz2}~\cite{Stipsicz}.

  The following natural case is the connected sum along embedded surfaces
  of genus $g=2$. Morgan and Szab\'o~\cite{Sz}~\cite{MSz} treated the case
  when the self-intersection of one of the surfaces is $1$ (and the other
  $-1$). The author~\cite{paper}~\cite{Thesis} 
  has solved the case in which the surfaces
  have self-intersection zero and are odd in homology, giving a number
  of nice applications. For this he used the Fukaya-Floer homology 
  as developed in~\cite{HFF}. 

  The case we are going to deal with in this paper is the connected sum
  along embedded surfaces of higher (arbitrary) genus $g$, self-intersection
  zero and odd, giving a positive answer to the conjectures in the 
  literature~\cite[conjecture 7.2]{MSz}~\cite[conjecture 5.27]{Thesis}. 

  Let $X$ be a smooth, compact, oriented four-manifold with $b^+ >1$
  and $b^+ -b_1$ odd. For 
  any $w \in H^2(X;\ZZ)$, $D^w_X$ will denote the corresponding Donaldson
  invariant~\cite{DK}~\cite{KM}, which is defined as a linear functional on 
  ${\AA(X)}= \Sym{*}(H_0(X)\oplus H_2(X)) \otimes \bigwedge^* H_1(X)$
  ($H_*(X)$ will always denote homology with rational coefficients, and
  similarly for $H^*(X)$). Let $x \in H_0(X)$ be the class of a point.
  Then Kronheimer and Mrowka~\cite{KM} define $X$ to be of simple type 
  (with respect to $w$) when $D_X^{w}((x^2 -4)z)=0$ for
  all $z \in {\AA(X)}$, and in that case define 
$$ 
  \DD_X^{w}(z)=D_X^{w}((1+{x \over 2})z),
$$
  for all $z \in \Sym{*}H_2(X)$. The series $\DD_X^{w}(e^{t\a})$, $\a \in
  H_2(X)$, is even or odd depending on whether
$$
  d_0=d_0(X,w)= -w^2- {3\over 2}(1-b_1+b^+)
$$ 
  is even or odd.
  When $b_1=0$ and $b^+>1$, $X$ is of simple
  type with respect to some $w$ if and only if it is so with respect to
  any $w$. In such case, $X$ is just called of simple type.

\begin{prop}[\cite{KM}]
\label{prop:KM}
  Let $X$ be a manifold of simple type with $b_1=0$ and 
  $b^+ >1$ and odd. Then we have 
  $$
     \DD_X^{w}(e^{\a})= e^{Q(\a) /2} \sum a_{j,w} \, e^{K_j \cdot \a}
  $$
  for finitely many $K_j \in H^2(X;\ZZ)$ (called basic
  classes) and $a_{j,w}=(-1)^{{K_j \cdot {w} +{w}^2}
     \over 2} a_j$, where $a_j$ are rational numbers (the collection is empty
  when the invariants all vanish). These classes are lifts
  to integral cohomology of ${w}_2(X)$. Moreover, for any embedded
  surface
  $\S \inc X$ of genus $g$ and with $\S^2 \geq 0$,
  one has $2g-2 \geq \S^2 +|K_j \cdot \S|$. 
\end{prop}

\begin{defn}[\cite{paper}]
\label{def:allowable}
  $(w,\S)$ is an {\bf allowable} pair if 
  $w, \S \in H^2(X; \ZZ)$, $w \cdot \S \equiv 1\pmod 2$ and $\S^2 =0$. 
  Then we define
  $$
    D^{(w,\S)}_X=D^w_X +D^{w+\S}_X.
  $$
  When $b^+=1$, the invariants depend on the metric through a 
  structure of walls and chambers~\cite{Kotschick}
  and therefore we have to specify the metric. As
  $w \cdot \S \equiv 1\pmod 2$, we have that
  $\S$ is contained in the closure of a unique chamber of
  type $(w,p_1)$ for every $p_1 <0$, $p_1 \equiv w^2 \pmod 4$.
  So, in the case that the invariants only depend on the 
  metric through the period point (for instance, simply-connected 
  manifolds~\cite{Kotschick2} 
  and $\S\x \CP^1$, with $\S$ a Riemann surface, which are
  all the cases we need for our arguments),
  we shall consider the invariants referring to
  the chambers defined by $\S$.
\end{defn}

  The series $D^{(w,\S)}_X(e^{t\a})$, $\a \in
  H_2(X)$, is even or odd according to whether 
  $d_0$ is even or odd. Since $(w+\S)^2 \equiv w^2 +2 \pmod 4$, 
  we can recover $D^w_X$ and $D^{w+\S}_X$ from $D^{(w,\S)}_X$. 

\begin{prop}[\cite{paper}]
\label{prop:DwS}
  Suppose $X$ is a manifold of simple type with $b_1=0$ and $b^+>1$
  and odd. Write the Donaldson series as $\DD_X^w (e^{\a})=
  e^{Q(\a)/2} \sum a_{j,w}\,
  e^{K_j \cdot \a}$. 
  Then setting $d_0=d_0(X,w)=-w^2-{3\over 2}(1+b^+)$ we have 
  $$
    D^{(w,\S)}_X(e^{\a})=e^{Q({\a})/2} \hspace{-8mm} \sum_{K_j \cdot
    \S \equiv 2\pmod 4} \hspace{-8mm} a_{j,w} e^{K_j \cdot \a} + 
    e^{-Q(\a)/2} \hspace{-8mm} \sum_{K_j
    \cdot \S \equiv 0\pmod 4} \hspace{-8mm} i^{-d_0} a_{j,w} e^{iK_j \cdot \a}
  $$
  So giving $\DD_X^{\,w}$ is equivalent to giving $D^{(w,\S)}_X$. 
\end{prop}

\begin{rem}
\label{rem:DwS}
 Under the conditions of proposition~\ref{prop:DwS}, we can prove in a similar
 fashion that
$$
 D^{(w,\S)}_X(x \,e^{\a}) = 2\, e^{Q({\a})/2} \hspace{-8mm} \sum_{K_j \cdot
    \S \equiv 2\pmod 4} \hspace{-8mm} a_{j,w} e^{K_j \cdot \a} -2 \,
    e^{-Q(\a)/2} \hspace{-8mm} \sum_{K_j
   \cdot \S \equiv 0\pmod 4} \hspace{-8mm} i^{-d_0} a_{j,w} e^{iK_j \cdot \a} 
$$
   and using that $ D^{(w,\S)}_X(x^a \,e^{\a+\S})= \sum {1 \over b!} 
   D^{(w,\S)}_X(\S^b x^a \,e^{\a})$, one has 
\begin{eqnarray*}
    D^{(w,\S)}_X(\S^b x^a \,e^{\a}) &= & e^{Q({\a})/2} \hspace{-8mm} \sum_{K_j 
    \cdot
    \S \equiv 2\pmod 4} \hspace{-8mm} a_{j,w} 2^a \,((\a+K_j)\cdot\S)^b \,
    e^{K_j \cdot \a} + \\
    & & + e^{-Q(\a)/2} \hspace{-8mm} \sum_{K_j
    \cdot \S \equiv 0\pmod 4} \hspace{-8mm} i^{-d_0} a_{j,w} (-2)^a \,
    ((- \a+ i \, K_j)\cdot\S)^b \, e^{iK_j \cdot \a}
\end{eqnarray*}
\end{rem}

\begin{defn}[\cite{paper}]
\label{def:permissible}
  We say that $(X,\S)$ is {\bf permissible} if $X$ is a smooth compact
  oriented four-manifold and $\S \inc X$ is an embedded Riemann surface
  of genus $g \geq 2$ and self-intersection zero such that $[\S] \in H_2(X;\ZZ)$
  is odd (its reduction modulo $2$ is non-zero, or equivalently, it is an
  odd multiple of a primitive homology class). So we can consider
  $w \in H^2(X;\ZZ)$ with $w \cdot \S \equiv 1 \pmod 2$. Then $(w,\S)$ 
  is an allowable pair. This implies that $b^+>0$. Let $N_{\S} \iso
  A =\S \x D^2$ be an open tubular neighbourhood of $\S$ and 
  set $\X=X - N_{\S}$. Then 
  $\bd \X=Y \iso \Y$ (but the isomorphism is not canonical). We consider
  one such isomorphism fixed and (when necessary) we furnish $\X$ with a
  cylindrical end, i.e. we consider $\X \cup (Y\x [0,\infty))$ (and keep
  on calling it $\X$).
\end{defn}

  We call {\bf identification} for $Y=\Y$ any (orientation preserving) 
  bundle automorphism $\p: Y \isom Y$. Up to isotopy, $\p$ depends only 
  on the isotopy class of the induced diffeomorphism on $\S$ and on an
  element of $H^1(\S;\ZZ)$.

\begin{defn}[\cite{paper}]
\label{def:sum}
  Let $(X_1,\S_1)$ and $(X_2,\S_2)$ be permissible, with $\S_1$ and 
  $\S_2$ of the same genus $g$. We pick orientations
  so that $\bd \X_1 = -\bd \X_2=Y$ (minus means reversed orientation).
  Then $X=X(\p) = \X_1 \cup_{\p} \X_2 = X_1 \#_{\S} X_2$ is  
  called the {\bf connected 
  sum along $\S$} of $(X_1,\S_1)$ and $(X_2,\S_2)$ (with identification
  $\p$). It is a compact,
  naturally oriented, smooth four-manifold with an embedding $Y \inc X$
  such that $X-Y= \X_1 \sqcup \X_2 \subset X_1 \sqcup X_2$.
\end{defn}

  The induced homology classes $[\S_1]$ and $[\S_2]$ coincide and
  are induced by an embedded $\S \inc X$. Then $(X,\S)$ is permissible.
  Choose $w_i \in H^2(X_i;\ZZ)$, $i=1,2$, and $w \in H^2(X;\ZZ)$ such 
  that $w_i \cdot \S_i \equiv 1 \pmod 2$, $w \cdot \S \equiv 1 \pmod 2$, 
  in a compatible way (i.e. the restricition of $w$ to $\X_i \subset X$ 
  coincides with the restriction of $w_i$ to $\X_i \subset X_i$).
  Always $w^2\equiv w_1^2+w_2^2 \pmod 2$. Changing $w$ by $w+\S$ if 
  necessary, we can always suppose $w^2\equiv w_1^2+w_2^2 \pmod 4$.
  In general, we shall call $w$ all of them, not making explicit
  to which manifold they refer. Note also that if $b_1(X_1)=b_1(X_2)=0$ then
  $b_1(X)=0$ and $b^+(X)>1$.

\begin{rem}
\label{rem:natural-identif}
  There is a case when there is a preferred identification. Suppose
  $X_1$ and
  $X_2$ are complex surfaces, $\S$ is a complex curve of genus $g$ and
  there are holomorphic
  embeddings $\S \inc X_i$, with image $\S_i$, $i=1,2$, such
  that $[\S_i]$ is odd and has self-intersection zero.
  Then the holomorphic normal bundle to $\S_i$ 
  gives a preferred isomorphism
  $\bd N_{\S_i} \isom \S_i \x \SS^1 \iso \S \x \SS^1$, and hence a 
  preferred identification $\bd \X_1 \isom -\bd \X_2$.
\end{rem}

\begin{rem}
\label{rem:homology-orientations}
  Let $X_1$ and $X_2$ be as in definition~\ref{def:sum}
  and $X=X_1 \#_{\S} X_2$. Then homology
  orientations for $X_1$ and $X_2$ induce conically a homology orientation
  for $X$. To see this, define $H_+^2(\X_i)$ to be a maximal definite
  positive subspace for the intersection pairing of $\X_i$ restricted to the
  image of $H_2(\X_i) \ar H_2(\X_i, \bd \X_i) \iso H^2(\X_i)$. Then it is
  esay to see that the exact sequence for the pair $(X_i,\X_i)$ and the
  Mayer-Vietoris sequence for $X=\X_1 \cup \X_2$ yield the
  following exact sequences
  $$ 0 \ar H^1(X_i) \ar H^1(\X_i) \ar H^2(A, \bd A) \ar H^2_+(X_i) 
     \ar H_+^2(\X_i) \ar 0 
  $$
  $$ 0 \ar  H^1(X) \ar H^1(\X_1) \oplus H^1(\X_2) \ar H^1(Y) \ar  
     H^2_+(X)  \ar H_+^2(\X_1) \oplus H^2_+(\X_2) \ar 0,
  $$
  where $A= \S \x D^2$ as before. We fix an orientation of $H^1(Y)$ (it
  does not change when we reverse the orientation of $Y$).
  So homology orientations for $X_1$ and $X_2$ induce a homology 
  orientation for $X$. 
\end{rem}

Let $\cH= \{  D \in H_2(X) /D|_Y=k[\SS^1]\in H_1(Y),
 \hbox{ some $k \in \QQ$} \}$. This subspace
of $H_2(X)$ contains the image of $H_2(\X_1) \oplus H_2(\X_2)$.
For every $D \in \cH$, choose $D_i \in H_2(X_i)$ agreeing with $D$ 
(i.e. $D_i|_{\X_i}=D|_{\X_i}$, $i=1,2$) 
and with $D^2 =D_1^2 +D_2^2$. Moreover, we
can suppose that the map $D \mapsto (D_1, D_2)$ is linear. Actually, once
chosen one of such maps, any other is of the form $D \mapsto
 (D_1+r\S, D_2 -r\S)$,
for a (rational) number $r$. 
Let us state now our main result.

\begin{thm}
\label{thm:main}
  Suppose $(X_1,\S_1)$ and $(X_2,\S_2)$ are permissible with $\S_1$ and 
  $\S_2$ of the same genus $g \geq 2$. Suppose also that $X_1$, $X_2$ 
  have both $b_1=0$ and $b^+>1$ and are of simple type. Let\/
  $\DD^w_{X_1}(e^{\a})=e^{Q(\a)/2}\sum a_{j,w} e^{K_j\cdot \a}$ and\/
  $\DD^w_{X_2}(e^{\a})=e^{Q(\a)/2}\sum b_{k,w} e^{L_k\cdot \a}$.
  Choose any identification $\p$ and let $X=X(\p)=X_1 \#_{\S} X_2$ 
  be the connected sum along $\S$, with the induded homology orientation.
  Suppose finally that $X$ is of simple type. Let $w \in H^2(X;\ZZ)$,
  $w_i \in H^2(X_i;\ZZ)$, $i=1,2$, in a compatible way, such 
  that $w_i \cdot \S_i \equiv 1 \pmod 2$, $w \cdot \S \equiv 1 \pmod 2$
  and $w^2\equiv w_1^2+w_2^2 \pmod 4$.    
  For every $D \in \cH$, we choose $D_i \in H_2(X_i)$ agreeing
  with $D$ satisfying $D^2 =D_1^2 +D_2^2$ in such a way that 
  the map $D \mapsto (D_1,D_2)$ is linear. Then
$$
   \DD^w_X(e^{tD}) = e^{Q(tD)/2}(\hspace{-5mm} \sum_{K_j \cdot \S=L_k \cdot
    \S = 2g-2}\hspace{-5mm}  -2^{7g-9} a_{j,w}b_{k,w} \, e^{(K_j \cdot D_1 + 
    L_k \cdot D_2 +2\S\cdot D)t} + 
$$
\begin{equation}
     + \hspace{-5mm}\sum_{K_j \cdot \S=L_k \cdot \S= -(2g-2)}
    \hspace{-5mm} (-1)^g \,2^{7g-9} a_{j,w}b_{k,w} \, e^{(K_j \cdot
    D_1 +  L_k \cdot D_2-2\S\cdot D)t}).
\label{eqn:main}
\end{equation}
\end{thm}

\begin{rem}
  If we do not assume $w^2 \equiv w_1^2 +w_2^2 \pmod 4$, we get an extra
  factor $\f=(-1)^{(g-1)(w^2 -w_1^2 -w_2^2)/2}$ in front of 
  formula~\eqref{eqn:main}.
\end{rem}

\begin{rem}
\label{rem:signs}
  The reason for the sign is easy to work out.
  First, $w^2 \equiv w_1^2+w_2^2 \pmod 2$. Also,
  $$
  -{3 \over 2}(1-b_1(X)+b^+(X))= -{3 \over
  2}(1-b_1(X_1)+b^+(X_1))-{3 \over 2}(1-b_1(X_2)+b^+(X_2)) -3(g-1).
  $$
  Therefore, $d_0(X,w) \equiv d_0(X_1,w)+d_0(X_2,w)+g-1 \pmod 2$. 
  Now the sign comes from the fact that
  the coefficient for the basic class $-\k$ is $(-1)^{d_0}c_{\k}$, being
  $c_{\k}$ the coefficient for the basic class $\k$.
\end{rem}

\begin{rem}
  If we are in the conditions of theorem~\ref{thm:main}, but $g=1$, we get
  a slightly different answer~\cite[chapter 4]{Thesis}~\cite{MSz}.
  For all basic classes, it is $K_j \cdot \S=L_j \cdot \S=0$, and
$$
   \DD^w_X(e^{tD}) = e^{Q(tD)/2}(\sum_{K_j,L_k} -{1 \over 4} a_{j,w}b_{k,w}
   \, e^{(K_j \cdot D_1 + L_k \cdot D_2 +2\S\cdot D)t} + 
$$
$$
 + \sum_{K_j,L_k} -{1 \over 4} a_{j,w}b_{k,w}
   \, e^{(K_j \cdot D_1 + L_k \cdot D_2 -2\S\cdot D)t} +
   \sum_{K_j,L_k} -{1 \over 2} a_{j,w}b_{k,w}
   \, e^{(K_j \cdot D_1 + L_k \cdot D_2)t}.
$$
\end{rem}

\begin{cor}
\label{cor:main}
  Suppose we are in the conditions of theorem~\ref{thm:main}. 
  Write\/ $\DD_X(e^{\a})=e^{Q(\a)/2} \sum c_{\k}
  e^{\k\cdot \a}$, 
  $\DD_{X_1}(e^{\a})=e^{Q(\a)/2}\sum a_j e^{K_j\cdot \a}$ and\/
  $\DD_{X_2}(e^{\a})=e^{Q(\a)/2}\sum b_k e^{L_k\cdot \a}$ for the Donaldson
  series for $X$, $X_1$ and $X_2$, respectively. 
  Then given any pair $(K,L) \in H^2(\X_1;\ZZ)
  \oplus H^2(\X_2;\ZZ)$, we have
  $$
   \sum_{\{\k/ \k|_{\X_1}=K, \, \k|_{\X_2}=L\}}\hspace{-4mm} c_{\k} = 
   (\pm 1)^{g-1}
    2^{7g-9}\,(\sum_{K_j|_{\X_1}=K} a_j) \cdot  (\sum_{L_k|_{\X_2}=L}
    b_k)
  $$
  whenever $K|_Y = L|_Y = \pm (2g-2) \PD [\SS^1]$. Otherwise, the left hand
  side is zero.
\end{cor}

\begin{pf}
  Allowing $D$ to vary in $\cH$, formula~\eqref{eqn:main} gives
  \begin{equation}
   \sum_{\{\k/ \k|_{\X_1}=K, \, \k|_{\X_2}=L\}}\hspace{-4mm} c_{\k,w} 
   e^{(\k\cdot D)t}=
   \sum_{{K_j|_{\X_1}=K} \atop {L_k|_{\X_2}=L}}
    - (\pm 1)^{g-1}
    2^{7g-9}\, a_{j,w} b_{k,w} \, e^{(K_j \cdot
    D_1 +  L_k \cdot D_2 \pm 2\S\cdot D)t}.
  \label{eqn:signs}\end{equation}
  We cannot have more precise information on $\k$ as we cannot evaluate
  $\DD_X^w$ on all $D \in H_2(X)$. Now take $w= \PD [D] \in H^2(X;\ZZ)$,
  $w_1 =\PD [D_1] \in H^2(X_1;\ZZ)$ and $w_2 =\PD [D_2] \in H^2(X_2;\ZZ)$
  with $w^2=w_1^2+w_2^2$ and $w \cdot \S =1$. Substitute 
  in~\eqref{eqn:signs} $t=\pi \, i/2$ and multiply by $(-1)^{w^2 / 2}=
  (-1)^{w_1^2/2} \, (-1)^{w_2^2 / 2}$, to get the sought expression.
\end{pf}

\begin{rem}
\label{rem:adm}
In theorem~\ref{thm:main} we cannot hope for having a similar formula for
classes $D$ such that $D|_Y$ is not a multiple of $[\SS^1]$ in $H_1(Y)$.
This is due to the fact that we cannot find $D_i \in H_2(X_i)$ agreeing
with $D$ (i.e. $D_i|_{\X_i}=D|_{\X_i}$, for $i=1,2$). We would need to
relate the invariant $\DD_X^w(e^{tD})$ with invariants of the
form $\DD_{\tilde X_i}^w(e^{tD_i})$, for suitable manifolds $\tilde
X_i$ containing $\X_i$, such that $D_i|_{\X_i}=D|_{\X_i}$. This was done
in~\cite[theorem 10]{paper} for the case $g=2$ (see also conjecture
in section~\ref{sec:conj}).

This limitation prevents us from having more general results. For example,
we do not know whether (under the conditions of
theorem~\ref{thm:main}) there are basic classes $\k$ for $X=X_1 \#_{\S}
X_2$ such that $|\k \cdot \S|<2g-2$ 
or not (compare~\cite[corollary 11]{paper}). 

In some cases, theorem~\ref{thm:main} is all that we need to find explicitly
the basic classes for $X$. This is due to the fact that there is a 
subspace $V \subset H_2(X)$ where all the basic classes vanish,
such that $H_2(X)=\cH \oplus V$. In~\cite[definition 4.1]{MSz}, 
Morgan and Szab\'o define admissible identification, which is a 
condition which implies that $X=X_1 \#_{\S} X_2$ is of simple type and
such $V$ exists. 
\end{rem}

\begin{cor}
  Suppose there exists a subspace $V \subset H_2(X)$ where all the basic 
  classes vanish such that $H_2(X)=\cH \oplus V$. Then there are no basic 
  classes $\k$ for $X$ such that $|\k \cdot \S |<2g-2$.
  The basic classes for $X$ are indexed by pairs of basic classes
  $(K_i, L_j)$ for $X_1$ and $X_2$ respectively, such that $K_i
  \cdot \S =L_j \cdot \S = \pm (2g-2)$.
\end{cor}

\begin{rem}
  Corollary~\ref{cor:main} agrees with the results of the kind for the 
  Seiberg-Witten invariants~\cite[section 7.3]{Thesis}.
  Morgan, Szab\'o and Taubes~\cite{MSzT} have proved the analogous
  result to corollary~\ref{cor:main} for the Seiberg-Witten basic
  classes (not the part corresponding to basic classes $\k$ for $X$
  with $|\k \cdot \S| <2g-2$). Both results are equivalent supposing
  true the conjecture of Witten~\cite{Witten} about the relationship
  of Donaldson and Seiberg-Witten invariants.
\end{rem}

Our last result is

\begin{thm}
\label{thm:finite}
  Let $S= \S \x \CP^1$, $w= \PD [\CP^1] \in H^2(S;\ZZ)$. Suppose that
  $S$ is of finite type of order $n \geq 1$ with respect to $w$ and
  $w+\S$ for the metrics defined by $\S$ (i.e. $\Dws_S((x^2-4)^n z)=0$, 
  for all $z \in \AA (S)$). Then for any  
  $(X,\S)$ permissible, $X$ is of finite type 
  of order at most $n$ with respect to any $w \in H^2(X;\ZZ)$ with 
  $w \cdot \S \equiv 1 \pmod 2$  (and for metrics defined by $\S$,
  in the case $b^+=1$). 
\end{thm}

\noindent {\em Acknowledgements:\/} Thanks to my D. Phil.\
supervisor Simon Donaldson, for many helpful ideas.
Also I am very grateful to the Mathematics 
Department in Universidad de M\'alaga for their hospitatility and support.

\section{Gluing theory}
\label{sec:gluing}

  Now we are going to develop the gluing theory necessary to prove 
  theorem~\ref{thm:main}. The set up is as follows, $X=\X_1 \cup_Y
  \X_2$, where $\bd \X_1=-\bd \X_2=Y$, an oriented three-manifold
  (later $Y=\Y$), $w \in H^2(X;\ZZ)$, $w_i=w|_{X_i} \in H^2(\X_i;\ZZ)$,
  and $D \in H_2(X)$ a homology class. So
  $D|_Y \in H_1(Y)$.
  We want to evaluate $D_X^w(D^d)$, the invariant being linear, we
  may multiply $D$ by any non-zero rational number, and hence suppose
  that $D|_Y$ is either zero or primitive in $H_1(Y;\ZZ)$.
  Now we represent $D$ by a cycle $D \subset X$ and put $\D_i =
  D \cap \X_i$, which we shall write formally as
  $D= \D_1+\D_2$, $\D_i \subset \X_i$. We can suppose
  $\bd \D_1=-\bd\D_2=\g$, with $\g \subset Y$ an embedded curve in $Y$,
  so when we give $\X_1$ a cylindrical end, we have
  $\D_1 \cap (Y\x [0,\infty)) = \g \x [0,\infty)$ (and analogously for
  $\X_2$).

\begin{prop}[\cite{HF}~\cite{HFF}~\cite{Thesis}]
\label{prop:Fukaya}
  Suppose $w|_Y$ odd. Then we have one of the following cases:
\begin{itemize}
  \item $D|_Y =0$  in $H_1(Y;\ZZ)$. Represent $D$ by a cycle so 
  $D= \D_1+\D_2$, $\D_i \subset \X_i$,
  $\bd \D_1=\bd\D_2=\emptyset$. Consider the Floer homology
  groups~\cite{HF}
  $HF_*(Y)$ (graded mod $4$). Then $(\X_i, \D_i)$ define relative 
  invariants
  $\p^{w_1}(\X_1,e^{t\D_1}) \in HF_*(Y) \otimes \QT$, 
  $\p^{w_2}(\X_2,e^{t\D_2})\in HF_*(-Y) \otimes \QT$. 
  There is a natural pairing
  $HF_*(Y) \otimes HF_*(-Y) \ar \QQ$, such that 
  \begin{equation}
     D_X^{(w,\S)}(e^{tD})=
     <\p^{w_1}(\X_1,e^{t\D_1}),\p^{w_2}(\X_2,e^{t\D_2})>. \label{eqn:sym}
  \end{equation}
  \item $D|_Y \neq 0$. Substitute $D$ by a rational
  multiple if necessary so that $D|_Y \in H_1(Y;\ZZ)$ is
  a primitive element in $H_1(Y;\ZZ)$.
  Represent $D$ by a cycle so $D= \D_1+\D_2$, $\D_i \subset \X_i$,
  $\bd \D_1=-\bd\D_2=\g$, with $\g \subset Y$ an embedded curve in $Y$. 
  Consider the Fukaya-Floer homology groups~\cite{HFF} 
  $HFF_*(Y,\g)$ 
  (graded mod $4$). Then $(\X_i, \D_i)$ define relative invariants
  $\p^{w_1}(\X_1,e^{t\D_1}) \in HFF_*(Y,\g)$, 
  $\p^{w_2}(\X_2,e^{t\D_2})\in HFF_*(-Y,-\g)$. There is a natural pairing
  $HFF_*(Y,\g) \otimes HFF_*(-Y,-\g) \ar \QT$, such that 
  \begin{equation}
     D_X^{(w,\S)}(e^{tD})=
     <\p^{w_1}(\X_1,e^{t\D_1}),\p^{w_2}(\X_2,e^{t\D_2})>. \label{eqn:sym2}
  \end{equation}
\end{itemize}
  When $b^+=1$, the invariants are calculated for a long neck, i.e. we
  refer to the invariants defined by $\S$.
\end{prop}

\begin{pf}
  \begin{itemize}
  \item As explained in~\cite{HF}, the Floer homology groups $HF_*(Y)$ are
  well-defined since $w|_Y$ is odd (this rules out problems with 
  flat reducible connections on $Y$). Also, there are invariants
  $\p^{w_1}(\X_1,(\D_1)^n) \in HF_*(Y)$, 
  $\p^{w_2}(\X_2,(\D_2)^m)\in HF_*(-Y)$ such that 
  $D_X^{(w,\S)}((\D_1)^n(\D_2)^m)= 
  <\p^{w_1}(\X_1,(\D_1)^n),\p^{w_2}(\X_2,(\D_2)^m)>$.
  Now we write 
  $$ \p^{w_1}(\X_1,e^{t\D_1})= \sum {t^n \over n!} \p^{w_1}(\X_1,(\D_1)^n),$$
  from where the statement of the theorem.
  \item Analogously, in~\cite{HFF} the Fukaya-Floer homology groups 
  $HFF_*(Y,\g)$ are defined when $w|_Y$ is odd. Associated to
  $(\X_i, \D_i)$, there are invariants
  $\p^{w_1}(\X_1,\D_1) \in HFF_*(Y,\g)$ and 
  $\p^{w_2}(\X_2,\D_2) \in HFF_*(-Y,-\g)$, where
  $\p^{w_1}(\X_1,\D_1)$ is represented by a Fukaya-Floer chain
  $(\p^{w_1}(\X_1,(\D_1)^n))_{n \geq 0}$, where $\p^{w_1}(\X_1,(\D_1)^n)
   \in HF_*(Y)$ (and analogously for $(\X_2,\D_2)$), such that
  $D_X^{(w,\S)}(D^r)= \sum_{0 \leq n \leq r} {r \choose n} 
  <\p^{w_1}(\X_1,(\D_1)^n),\p^{w_2}(\X_2,(\D_2)^{r-n})>$.
  So we write formally
  $$ \p^{w_1}(\X_1,e^{t\D_1})= \sum {t^n \over n!} \p^{w_1}(\X_1,(\D_1)^n).$$
  \end{itemize}
\end{pf}

Now we particularize to the case which concerns us, $Y=\Y$. Conjugation in
the second factor produces
an isomorphism $Y \iso (-Y)$ (also $(Y,\g) \iso (-Y, -\g)$).
As explained in~\cite{HFF}~\cite{Thesis}, 
$HFF_*(Y,\g)$ is
the limit of a spectral sequence whose $E_3$-term is $HF_*(Y) \otimes
\hH_*(\CP^{\infty})$ (here $\hH_*(\CP^{\infty})$  
means the natural completion of 
$H_*(\CP^{\infty})$, i.e. $\QT$), with differencital $d_3$ given by 
$$
   \mu(\g):HF_i(Y) \otimes H_j(\CP^{\infty}) \ar HF_{i-3}(Y) \otimes
   H_{j+2}(\CP^{\infty}).
$$

Let $\g=\pt \x \SS^1 \subset \Y$. Now all the
differentials in the $E_3$ term of the spectral sequence are of the
form $H_{\mathrm{odd}}(\M) \ar H_{\mathrm{even}}(\M)$ and
$H_{\mathrm{even}}(\M) \ar H_{\mathrm{odd}}(\M)$. The boundary
cycle $\g=\SS^1$ is invariant under the action of the group
$\Diff (\S)$ on $Y=\Y$, so the differentials commute with the
action of $\Diff (\S)$. As there are elements $\rho \in
\Diff (\S)$ acting as $-1$ on $H^1(\S)$, we have that $\rho$
acts as $-1$ on $H_{\mathrm{odd}}(\M)$ and as $1$ on
$H_{\mathrm{even}}(\M)$. Therefore the differentials are zero and the
spectral sequence degenerates in the third term. 
This implies that 
$HFF_*(Y,\g) = HF_*(Y) \otimes \hH^*(\CP^{\infty}) =V [[t]]$, 
where $V=HF_*(Y)$.
Now the relative invariants for $(\X_1, \D_1)$ can be written as 
$$ 
  \p^{w_1}(\X_1,e^{t\D_1})= \sum {t^n \over n!} 
  \p^{w_1}(\X_1,(\D_1)^n) \in V[[t]]
$$
where $\p^{w_1}(\X_1,(\D_1)^n) \in HF_*(Y)$ has perfect meaning.
Under the isomorphism $HFF_*=
HFF_*(Y, \g)=V[[t]]$, $HFF_*$ becomes a $\QT$-module
and
$HFF_* \otimes HFF_* \ar \QT$ is $\QT$-bilinear.

\begin{cor}
\label{cor:concl}
\begin{enumerate}
  \item There is a (rational) vector space $V=HF_*(Y)$ endowed with a
    bilinear form such that for every permissible $(X,\S)$,
    $w \in H^2(\X;\ZZ)$ with $w|_Y$ an odd multiple of $\PD [\SS^1]$,
    and cycle 
    $\D \subset \X$ with $\bd \D =\emptyset$, we have
    $\p^w(\X,e^{t\D}) \in V[[t]]$. For $X=\X_1 \cup_Y \X_2$, 
    $D=\D_1 +\D_2$, $\bd \D_1 =\bd \D_2 =\emptyset$, we have
    $$ D_X^{(w,\S)}(e^{tD})=
     <\p^{w_1}(\X_1,e^{t\D_1}),\p^{w_2}(\X_2,e^{t\D_2})>. $$
  \item There is a canonical isomorphism $HFF_*(Y,\SS^1)\iso V [[t]]$, 
    such that for every permissible $(X,\S)$,
    $w \in H^2(\X;\ZZ)$ with $w|_Y$ an odd multiple of $\PD [\SS^1]$,
    and cycle $\D \subset \X$ with $\bd \D =\SS^1$, we have
    $\p^w(\X,e^{t\D}) \in V[[t]]$. For $X=\X_1 \cup_Y \X_2$, 
    $D=\D_1 +\D_2$, $\bd \D_1 =-\bd \D_2 =\SS^1$, we have
    $$ D_X^{(w,\S)}(e^{tD})=
     <\p^{w_1}(\X_1,e^{t\D_1}),\p^{w_2}(\X_2,e^{t\D_2})>. $$
\end{enumerate}
\end{cor}

\begin{prop}
\label{prop:H(M)}
  Let $\M$ be the moduli space of odd degree rank two 
  stable vector bundles on $\S$, which is a smooth variety~\cite{Alaistar}.
  Then there is an isomorphism 
  $$ 
   HF_*(Y) \iso  H_*(\M)
  $$
  as vector spaces (we are using rational coefficients), where
  we reduce the grading of $H_*(\M)$ modulo $4$.
\end{prop}

\begin{pf}
Dostoglou and Salamon~\cite{DS} prove
$HF_*(\Y) \iso HF_*^{\mathrm{symp}}(\M)$. It is the particular case where
we consider $\p=\hbox{id}:\S \ar \S$, in which the mapping torus of $\p$ is
$\Y$. As explained in the introduction of~\cite{DS}, $\M$ is connected,
simply connected and $\pi_2(\M)=\ZZ$, so the groups
$HF_*^{\mathrm{symp}}(\M)$ are well-defined.
Now $HF_*^{\mathrm{symp}}(\M) \iso H_*(\M)$ is a standard result
obtained by Floer himself~\cite{Floer} for proving the Arnold conjecture.
\end{pf}

There is a map $$\mu: H_*(\S) \to H^{4-*}(\M)$$
given by $\mu(\a)= -{1 \over 4} \,
p_1(\cV) / 4$, where $\cV \to \S \x \M$ is the associated
universal $SO(3)$-bundle, $p_1(\cV) \in H^4(\S\x\M)$ its first
Pontrjagin class.
Fix a basis $\{\g_i\}$ of $H_1(\S)$. Let $a=\mu(\S)$, $b=\mu(x)$, 
$c_i=\mu(\g_i)$. These elements generate $H^*(\M)$ as a 
ring~\cite{Alaistar}~\cite{Thaddeus}.
So there is a basis for $V=H^*(\M)$ with elements of the form
$$
   f_{\a}=a^nb^mc_{i_1}\cdots c_{i_r} \in V,
$$
for a finite set of indices of the
form $\a=(n,m; i_1,\ldots,i_r)$, $n,m \geq 0$, $r \geq 0$, $1 \leq
i_1 < \cdots < i_r \leq 2g$. Let $N=\dim V$. We order the set of
indices $\{\a\}$ so we identify such set with $\{ 1, \ldots, N\}$
and write $1 \leq \a \leq N$ in general.

Let $I$ be the ideal of $H^*(\M)$ 
generated by $c_1, \ldots, c_{2g}$. Then the
elements $a^n b^m$, $0\leq n,m <g$, generate the quotient $H^*(\M)/I$
(see~\cite{Alaistar}). So we
can suppose these elements are the first $g^2$ 
elements in the basis $\{f_{\a}\}$, i.e. for $1 \leq \a \leq g^2$.

The intersection pairing in $H^*(\M)$ is given by 
\begin{equation} 
  <f_{\a}, f_{\b}> = <f_{\a}\cup f_{\b}, [\M]>.
\label{eqn:fa}\end{equation}
Therefore the intersection matrix $( <f_{\a}, f_{\b}>  )$ is invertible.

Here we recall that we have defined the manifold $A=\S \x D^2$,
with boundary $Y=\Y$, and let $\De= \pt \x D^2 \subset A$ be the 
horizontal slice with $\bd \De=\SS^1$. Put $w=\PD [\De] 
\in H^2(A;\ZZ)$. Put $\AA(\S)= \Sym{*}(H_0(\S)
\oplus H_2(\S)) \otimes \bigwedge^* H_1(\S)$, so there is
a natural map $\AA(\S) \ar \AA(X)$, whenever $\S \inc X$.
Then we define
\begin{eqnarray*}
   z_{\a}&=& \S^n x^m \g_{i_1} \cdots \g_{i_r} \in \AA(\S), \\
   e_{\a}&=& \p^w (A, z_{\a} e^{t\De}) \in HFF_*(Y, \SS^1)= V[[t]].
\end{eqnarray*}

\begin{lem}
\label{lem:ea}
  The intersection matrix $( <e_{\a}, e_{\b}>  )$ (with coefficients in
  $\QT$) is invertible.
  Therefore, $\{e_{\a}\}$ is a basis for $HFF_*(Y)=V[[t]]$.
\end{lem}

\begin{pf}
  As the elements $f_{\a} \in H^*(\M)$ have an integer degree between $0$
  and $6g-6=\dim \M$, we can reorder the basis 
  $\{f_{\a}\}$ such that the degree goes increasing (we use this
  special ordering only in this lemma). Now 
  $z_{\a}= \S^n x^m \g_{i_1} \cdots \g_{i_r}$, $z_{\b}= \S^{n'} x^{m'} 
  \g_{i'_1} \cdots \g_{i'_{r´}}$ and 
  $$
    <e_{\a}, e_{\b}> =  <\p^w (A, z_{\a} e^{t\De}),
    \p^w (A, z_{\b} e^{t\De})>= \Dws_{\S \x \CP^1} (z_{\a}z_{\b}e^{t
    \,\CP^1}),
  $$
  where $w=\PD [\CP^1] \in H^2(\S\x\CP^1;\ZZ)$.
  The matrix $( <e_{\a}, e_{\b}>  )$  has coefficients in $\QT$, so it
  is invertible if and only if its determinant is a unit in $\QT$, i.e.
  when we put $t=0$ we obtain an invertible matrix with rational
  coefficients. Now
  $$
    <e_{\a}, e_{\b}>|_{t=0} =   \Dws_{\S \x \CP^1} (z_{\a}z_{\b}).
  $$
  The lowest dimension of the moduli spaces of anti-self-dual connections
  for $\S \x \CP^1$
  is $6g-6$, so if  $\deg z_{\a}+ \deg z_{\b} < 6g-6$ then 
  $<e_{\a}, e_{\b}>|_{t=0} =    0$. The moduli space of dimension $6g-6$ 
  corresponds to $p_1=0$ and $w=\PD[\CP^1]$. 
  All of these connections are flat and
  irreducible, actually pull-back of flat connections on $\S$ with $w|_{\S}$
  odd, so the corresponding moduli space is isomorphic to $\M$.
  Thus if $\deg z_{\a}+ \deg z_{\b} = 6g-6$ then 
  $<e_{\a}, e_{\b}>|_{t=0} =    <f_{\a}, f_{\b}>$.
  The matrix $ (<f_{\a}, f_{\b}>)$ is of the form
  $$
    \left( \begin{array}{ccccc}
    0 & 0 & \cdots & 0& A_0 \\
    0  & 0 & \cdots & A_1 & 0 \\
    \vdots & \vdots & \ddots & \vdots & \vdots \\
    A_{6g-6} & 0& \cdots &0 &0  
    \end{array} \right)
  $$
  where $A_i$ are submatrices corresponding to the intersection
  product $$H^i(\M) \otimes H^{6g-6-i}(\M) \ar \QQ.$$ So all $A_i$
  have non-vanishing determinant $\det A_i \in \QQ$. 
  Finally, we have that the matrix $ (<e_{\a}, e_{\b}>|_{t=0})$ is 
  $$
    \left( \begin{array}{ccccc}
    0 & 0 & \cdots & 0& A_0 \\
    0  & 0 & \cdots & A_1 & \hbox{*} \\
    \vdots & \vdots & \ddots & \vdots & \vdots \\
    A_{6g-6} & \hbox{*} & \cdots & \hbox{*} & \hbox{*}  
    \end{array} \right)
  $$
  and it is invertible.
\end{pf}

\section{Proof of Theorem~\protect\ref{thm:main}}

  By the above lemma, $\{ e_{\a}\}$ is a basis of $V[[t]]$, 
  so there is an isomorphism
  \begin{eqnarray*}
     V[[t]] &\ar& \QQ^N[[t]] \\
     \p    &\mapsto& (<\p, e_{\a}>)_{\a}.
  \end{eqnarray*} 
  The important feature is that if $(X_1,\S_1)$ is permissible, $w \in H^2(X_1
  ; \ZZ)$ with $w \cdot \S \equiv 1 \pmod 2$, $\D_1 \subset \X_1$ with
  $\bd \D_1 =\SS^1$, $D_1=\D_1 + \De$, then
  $\p=\p^w(\X_1, e^{t\D_1})$ goes to $(c_{X_1,\a}(t))_{1 \leq \a \leq N} 
  \in \QQ^N[[t]]$, where
  $$
     c_{X_1,\a}(t)=\Dws_{X_1} (z_{\a} e^{tD_1}).
  $$

The pairing in $V[[t]]$ corresponds through the isomorphism to a pairing
in $\QQ^N[[t]]$, which is $\QT$-bilinear, hence given by a matrix 
of $\QT$-coefficients $(M_{\a\b}(t))_{1 \leq \a,\b \leq N}$. This 
matrix is universal (only dependent on the data necessary for the 
construction of the Fukaya-Floer groups, i.e. $(Y,\g)$, and on the
chosen basis). 

Now if $(X_1, \S_1)$ and $(X_2, \S_2)$ are permissible,
let $X=X(\p)= X_1 \#_{\S} X_2$ (with an identification $\p$)
be the connected sum along $\S$, with the induded homology orientation.
Let $D \in H_2(X)$ with $D=\D_1 +\D_2$, $\bd \D_1= - \bd \D_2= \SS^1$. 
Put $D_i=\D_i+\De$. Then (here $c_{X_1,\a}(t)=\Dws_{X_1} (z_{\a} e^{tD_1})$,
$c_{X_2,\b}(t)=\Dws_{X_2} (z_{\b} e^{tD_2})$),
\begin{equation}
  \Dws_X(e^{tD})=<\p^w(\X_1, e^{t\D_1}),\p^w(\X_2, e^{t\D_2})>=
  \sum_{1 \leq \a,\b \leq N} c_{X_1,\a}(t) M_{\a\b}(t) c_{X_2,\b}(t)
\label{eqn:f1} \end{equation}

Now suppose that $X_1$ has $b_1=0$. 
Then $c_{X_1,\a}(t)=0$ whenever $r>0$ (recall $\a=(n,m; i_1, \ldots, i_r)$). 
So the only non-zero coordinates
correspond to $z_{\a}=\S^n x^m$, $0 \leq n,m \leq g-1$.
Suppose furthermore $X_1$ of simple type with $b_1=0$ and $b^+>1$, so
for $z = (x^2-4)\S^n x^{m-2}$, 
$0 \leq n \leq g-1$, $2 \leq m \leq g-1$, $\Dws_{X_1}(z e^{tD_1})=0$,
hence changing the 
basis $\{z_{\a} \}$, all coordinates $c_{X_1,\a}(t)$ 
are zero except for the first $2g$ of them, corresponding
to $z_{\a}=\S^n$ and $z_{\a}=\S^n x$, $0 \leq n \leq g-1$.

\begin{lem}
\label{lem:z=0}
  Let $(X_1,\S_1)$ be permissible, with $\S_1$ of genus $g$, $X_1$
  of simple type with $b_1=0$ and $b^+>1$, $(w,\S)$ allowable and
  $D_1 \in H_2(X_1)$ with $D_1 \cdot \S =1$. Then $\Dws_{X_1}(
  e^{tD_1}z)=0$ for 
$$
  z=\left\{ \begin{array}{ll} 
   (1-{x \over 2})(\S+1)((\S+1)^2+4^2)((\S+1)^2+8^2) \cdots
   ((\S+1)^2+(2g-4)^2)  &  \hbox{$g$ even} \\
   (1+{x \over 2})(\S+1)(\S-3)(\S+5) \cdots
   (\S-(2g-3))  &  \hbox{$g$ odd}
  \end{array} \right.
$$
\end{lem}

\begin{pf}
  Suppose, for instance, $g$ even. By remark~\ref{rem:DwS}, for any
  polynomial $p(\S)$ in $\S$, 
  $$
    \Dws_{X_1}((1-{x \over 2}) p(\S)e^{tD_1})=
    2 e^{-Q({tD_1})/2} \hspace{-8mm} \sum_{K_j \cdot
    \S \equiv 0 \pmod 4} \hspace{-8mm} i^{-d_0} p((-D_1 +i\, K_j) 
    \cdot \S) a_{j,w} e^{itK_j \cdot D_1}.
  $$
  Now $D_1 \cdot \S=1$, so the expression above vanishes when $p(\S)$ 
  has roots $-1, -1 \pm 4i , \ldots -1 \pm (2g-4)i$.
\end{pf}

Note that $z=(1 \pm {x \over 2}) p(\S)$ with $p$ of degree $g-1$.
Let us choose a basis $\{z_{\a}\}$ with $z_1, \ldots z_{2g-1}$ being the
elements 
$$
  1, \S, \S^2, \ldots, \S^{g-1}, x, \S x, \ldots, \S^{g-2} x,
$$
$z_{2g}=z$, also $z_{2g+1}, \ldots, z_{g^2}$ being the elements $(x^2-4)
\S^n x^{m-2}$, $0 \leq n \leq g-1$, $2 \leq m \leq g-1$ and
$e_{g^2+1}, \ldots, e_N$ having all $r>0$.
So, when $X_1$ is of simple type with $b_1=0$ and $b^+>1$, 
$$
   \p^w(\X_1, e^{t\D_1})\in \QQ^{2g-1}[[t]] \subset V[[t]],
$$
where $\QQ^{2g-1}$ is the orthogonal complement to $<e_{2g}, \ldots, e_N>$
in $V$.

Formula~\eqref{eqn:f1} reduces to (when both $X_i$ are of simple
type with $b_1=0$ and $b^+>1$)
\begin{equation}
  \Dws_X(e^{tD})=\sum_{1 \leq \a,\b\leq 2g-1} c_{X_1,\a}(t) M_{\a\b}(t) 
   c_{X_2,\b}(t),
\label{eqn:f2} \end{equation}
where
$$
  c_{X_1,\a}(t)=\left\{
    \begin{array}{ll} 
    \Dws_{X_1}(\S^n e^{tD_1}) \qquad & \hbox{if $0 \leq n \leq g-1$,
	     $\a=(n,0;)=n+1$} \\
    \Dws_{X_1}(\S^n x \,e^{tD_1}) \qquad & \hbox{if $0 \leq n \leq g-2$,
	     $\a=(n,1;)=n+g+1$} 
    \end{array} \right.
$$
$$
  =\left\{
    \begin{array}{l} 
    e^{Q(tD_1)/2} \hspace{-8mm} \sum\limits_{K_j \cdot
      \S \equiv 2\pmod 4} \hspace{-8mm} a_{j,w} (1+ K_j \cdot\S)^n
      e^{t K_j \cdot D_1} + e^{-Q(t D_1)/2} \hspace{-8mm} \sum\limits_{K_j
      \cdot \S \equiv 0\pmod 4} \hspace{-8mm} i^{-d_0} a_{j,w} (- 1+ i \, 
      K_j \cdot \S)^n e^{ti\,K_j \cdot D_1} \\
    2 e^{Q(tD_1)/2} \hspace{-9mm} \sum\limits_{K_j \cdot
      \S \equiv 2\pmod 4} \hspace{-9mm} a_{j,w} (1+ K_j \cdot\S)^n
      e^{tK_j \cdot D_1} -2 e^{-Q(tD_1)/2} \hspace{-9mm} \sum\limits_{K_j
      \cdot \S \equiv 0\pmod 4} \hspace{-9mm} i^{-d_0} a_{j,w} (- 1+ i \, 
      K_j \cdot \S)^n e^{tiK_j \cdot D_1} 
    \end{array} \right.
$$
and analogously for $c_{X_2,\b}(t)$.

So it is easy to find another basis $\{z_1, \ldots, z_{2g-1} \}$ (which
we do not write explicitly) 
spanning $\QQ^{2g-1}$ such that
$$
  c_{X_1,\a}(t)=\left\{
    \begin{array}{ll} 
    e^{Q(tD_1)/2} \sum\limits_{K_j \cdot
      \S =2p} a_{j,w} e^{tK_j \cdot D_1} \qquad
       & \hbox{if $p$ is odd} \\
    e^{-Q(tD_1)/2} \sum\limits_{K_j \cdot
      \S =2p} i^{-d_0}a_{j,w} e^{ti\, K_j \cdot D_1} \qquad
      & \hbox{if $p$ is even} 
    \end{array} \right.
$$
where $\a=1, \ldots, 2g-1$ corresponds to $p=g-1, -(g-1), g-2, -(g-2), \ldots,
0$. Formula~\eqref{eqn:f2} yields
\begin{equation}
  \Dws_X(e^{tD})=\sum_{1 \leq \a,\b\leq 2g-1} c_{X_1,\a}(t) M_{\a\b}(t) 
  c_{X_2,\b}(t) 
\label{eqn:f3} \end{equation}
with $c_{X_2,\b}(t)$ defined analogously to $c_{X_1,\a}(t)$, with 
the letter $q$ in the place of $p$. 

\begin{lem}
\label{lem:dia}
  The matrix $(M_{\a\b}(t))$ is diagonal.
\end{lem}

\begin{pf}
First we note that
$$
  \Dws_X(e^{tD +s\S})= \sum_{p \; odd} c_{X,\a}(t) e^{2ps} +
  \sum_{p \; even} c_{X,\a}(t) e^{2pi s}.
$$
We use equation~\eqref{eqn:f3} for $D+ {s \over t}\S$, $D_1 +
{s \over t}\S$ and $D_2$, so   
$$
  \Dws_X(e^{tD +s\S})= \hspace{-2mm}
  \sum_{1 \leq \a,\b\leq 2g-1 \atop p \; odd} \hspace{-2mm}
  c_{X_1,\a}(t) M_{\a\b}(t)  c_{X_2,\b}(t) e^{2ps}+ \hspace{-2mm}
  \sum_{1 \leq \a,\b\leq 2g-1 \atop p \; odd}  \hspace{-2mm}
  c_{X_1,\a}(t) M_{\a\b}(t)  c_{X_2,\b}(t) e^{2pi s}.
$$
Then 
$$
  c_{X\a}(t)=\sum_{\b} c_{X_1,\a}(t) M_{\a\b}(t)  c_{X_2,\b}(t).
$$
Let us see that $M_{\a\b}(t)=0$ unless $\b=\a$. Suppose, for instance,
that $p$ is odd.
If we write now $D=(\D_1+r\S)+ (\D_2-r\S)$, the left hand side of the 
expression above remains unchanged, but the right hand side is a sum of 
exponentials $e^{(2p-2q)r}$, $q$ odd, and $e^{(2p-2qi)r}$, $q$ even. So 
it when $q \neq p$, it is $M_{\a\b}(t)=0$.
\end{pf}

Formula~\eqref{eqn:f3} gives
$$
  \Dws_X(e^{tD})=\hspace{-8mm} 
    \sum_{K_j \cdot \S =L_k \cdot \S =2p,\; p \; odd}  \hspace{-8mm}  
     e^{Q(tD_1)/2+Q(tD_2)/2} M_{\a\a}(t) a_{j,w} b_{k,w}
     e^{tK_j \cdot D_1+tL_k \cdot D_2} +
$$
$$
    + \sum_{K_j \cdot \S =L_k \cdot \S=2p,\; p \; even} \hspace{-8mm}    
     e^{-Q(tD_1)/2-Q(tD_2)/2}  M_{\a\a}(t) i^{-d_0(X_1,w_1)-d_0(X_2,w_2)}
     a_{j,w} b_{k,w}
     e^{ti\,K_j \cdot D_1+ti\,L_k \cdot D_2}. 
$$

Obviously $D^2=D_1^2+D_2^2$. We are assuming that $X$ is of simple type 
and recall that $b_1(X)=0$ and $b^+(X)>1$. By remark~\ref{rem:signs}, 
$d_0(X,w)-d_0(X_1,w_1) - d_0(X_2,w_2)=w^2-w_1^2-w_2^2-3(g-1)\equiv
g-1 \pmod 4$ (since we are assuming $w^2 \equiv w_1^2+w_2^2 \pmod 4$), so
$$
  \DD^w_X(e^{tD})=e^{Q(tD)/2} \hspace{-2mm} 
    \sum_{K_j \cdot \S =L_k \cdot \S=2p \atop -(g-1) \leq p
    \leq g-1}  \hspace{-2mm}
      i^s M_{\a\a}(t) a_{j,w} b_{k,w}
     e^{tK_j \cdot D_1+tL_k \cdot D_2}
$$
where $s=0$ when $p$ is odd and $s=g-1$ when $p$ is even. 
As $g$ is not dependent on the particular manifolds, 
we absorb this factor into the matrix without affecting its universality.

This expression is valid for any $D \in H_2(X)$ such that $D|_Y=\SS^1$. 
We note that it does not change if we change $(D_1, D_2)$ for
$(D_1+r\S, D_2-r\S)$, as expected. This means that we only need to assume
the conditions: 
$D$ and $D_i$ coincide in $\X_i$, $i=1,2$, and $D^2=D_1^2+D_2^2$.
Now take a linear map $D \mapsto (D_1,D_2)$ from the subspace
$\cH=\{D \in H_2(X) / D|_Y=k [\SS^1], \hbox{ some $k$}\}$, satisfying the
former conditions. As the set of
$D \in H_2(X)$ with $D|_Y$ a non-zero multiple of $\SS^1$ is
dense in $\cH$, we have that for any $D \in \cH$,
\begin{equation}
  \DD^w_X(e^{tD})=e^{Q(tD)/2} \hspace{-2mm} 
    \sum_{K_j \cdot \S =L_k \cdot \S=2p}  \hspace{-2mm}
      M_{\a\a}( t(D\cdot \S) ) a_{j,w} b_{k,w}
     e^{tK_j \cdot D_1+t L_k \cdot D_2}
\label{eqn:f4} 
\end{equation}

\begin{lem}
\label{lem:dia2}
  $M_{\a\a}(t)=0$, except for $\a=1,2$.
\end{lem}

\begin{pf}
  Let $Y$ be the $K3$ surface, which is a manifold of simple type  
  with $b_1=0$ and $b^+>1$ (see~\cite{KM}). Consider 
  a tight embedded Riemann surface $\S' \inc Y$ of genus $g'<g$ and
  self-intersection zero. By definition of tightness (see~\cite{KM}),
  $(\S')^2=2g'-2$. To construct it, we consider an elliptic fibration for
  the $K3$ surface. Let $T$ be a generic fibre (which is a torus of
  self-intersection zero) and let $S$ be a section represented by a sphere
  of self-intersection $-2$. Then consider $S$ together with $g'$ generic
  fibres and smooth out the intersection points. Call the resulting Riemann 
  surface $\S'$. It is homologous to $S+g'T$, it 
  has genus $g'$ and self-intersection $2g'-2$, as required.

  We blow-up $Y$ at $2g'-2$ points in $\S'$, to get $X_1=Y \#(2g'-2)
  \overline{\CP}^2$. The proper transform of
  $\S'$ is a Riemann surface of genus $g'$ and self-intersection zero.
  Perform an internal connected sum with $g-g'$ homologically trivial tori
  to obtain a Riemann surface $\S_1$ of genus $g$ and self-intersection zero.
  If $E_1, \ldots, E_{2g'-2}$ are the exceptional divisors, $\S_1$ is
  homologous to $\S' +E_1 +\cdots +E_{2g'-2}$.

  Moreover, the basic classes of $X_1$ are $ \pm E_1 +\cdots  \pm
  E_{2g'-2}$. They all satisfy $\k \cdot \S_1 \leq 
  2g'-2$ and there is exactly one, $K =E_1 +\cdots +E_{2g'-2}$,
  satisfying the equality.

  Let $(X_2, \S_2)=(X_1, \S_1)$, and consider $X = X_1 \#_{\S}
  X_2$ with the preferred 
  identification of remark~\ref{rem:natural-identif} 
  (double of $X_1$ along $\S_1$). Then 
  $X$  splits off a ${\Bbb S}^2 \x {\Bbb S}^2$, so its invariants are
  zero (see~\cite[section 4.3]{MSzT} for this well-known phenomenon).
  As in the proof of lemma~\ref{lem:dia},
  $c_{X,\a}(t)=c_{X_1,\a}(t) M_{\a\a}(t)  c_{X_2,\a}(t)$.
  We proceed by induction from $g'=1, 2, \ldots, g-1$ and get
  $M_{\a\a}(t)=0$, for $\a \geq 3$.
\end{pf}

Formula~\eqref{eqn:f4} becomes
$$  
   \DD^w_X(e^{tD})=e^{Q(tD)/2} \hspace{-8mm} 
    \sum_{K_j \cdot \S =L_k \cdot \S=2g-2} \hspace{-8mm}   
      M_{11}( t(D\cdot \S) ) a_{j,w} b_{k,w}
     e^{tK_j \cdot D_1+t L_k \cdot D_2} +
$$
\begin{equation}
   + e^{Q(tD)/2} \hspace{-8mm} 
    \sum_{K_j \cdot \S =L_k \cdot \S=-(2g-2)} \hspace{-8mm}    
      M_{22}( t(D\cdot \S) ) a_{j,w} b_{j,w}
     e^{tK_j \cdot D_1+t L_k \cdot D_2}
\label{eqn:f5} \end{equation}

Now let us compute $M_{11}(t)$ and $M_{22}(t)$. We can do that with a 
particular example of connected sum
along a surface of genus $g$ in which the invariants are known for $X_1$,
$X_2$ and $X=X_1 \#_{\S} X_2$.

\begin{defn}
\label{def:Bg,Cg}
  We are going to define the following smooth manifolds:
\begin{itemize}
  \item  Let $S_n$ be the minimal
	 elliptic surface with no multiple fibres and 
	 geometric genus $p_g=n-1$. This is unique up to 
	 diffeomorphism~\cite{Friedman}. $S_1$ is the rational elliptic
	 surface, i.e. $S_1= \CP^2 \#Ê9 \overline{\CP}^2$.
  \item Let $\tilde S_1$ be the blow-up of $S_1$ at one point. Therefore,
   $\tilde S_1=\CP^2 \#Ê10 \overline{\CP}^2$. Consider a particular elliptic
   fibration for $S_1$ with section $\s$ 
   (of self-intersection $-1$) and fibre $F$. Call $E$ the
   exceptional divisor, so $\tilde S_1$ has an elliptic fibration
   with fibre the torus $T_1=F$ and there is another embedded torus $T_2$
   homologous to $\s+T_1-E$ with $T_2^2=0$ and $T_1 \cdot T_2 =1$.
  \item Let $B_g=\underbrace{\tilde S_1 \#_{T_1}\tilde S_1 \#_{T_1} \cdots
    \#_{T_1} \tilde S_1}_{g}$ 
   (connected sums along $T_1$'s with the preferred identification), 
   which is diffeomorphic to $S_g \# g \overline{\CP}^2$. It
   contains an embedded torus $T_1$ of self-intersection zero
   and a Riemann surface of genus $g$ (and
   self-intersection zero) made up gluing smoothly the $T_2$'s
   coming from each $\tilde S_1$. Actually, the elliptic surface $S_g$ has a 
   section $\s_g$ with $\s_g^2=-g$, and $\S_g$ can be taken to be the
   proper transform of $\s_g$. Clearly,
   $\S_g \cdot T_1 =1$, so $(B_g, \S_g)$ is permissible. $B_g$ is of simple
   type with $b_1=0$, $b^+>1$ (as $g \geq 2$).
  \item Let $C_g= B_g \#_{\S_g} B_g $ with the preferred identification.
   It contains a Riemann surface $\hat \S_2$ 
   of genus $2$ and self-intersection zero
   made up from gluing smoothly the $T_1$'s.
   If we perform instead the connected sum of two $\tilde S_1$ along $T_2$,
   we get $\hat B_2=\tilde S_1
   \#_{T_2} \tilde S_1 $ with an embedded Riemann surface $\hat \S_2$ 
   of genus $2$ and self-intersection zero, coming
   from smoothly gluing the $T_1$'s. Clearly $(\hat B_2,\hat \S_2) \iso 
   (B_2,\S_2)$. Now
  $$ 
     C_g= B_g \#_{\S_g} B_g=\underbrace{\hat B_2 \#_{\hat \S_2}  \cdots
   \#_{\hat\S_2} \hat B_2}_{g}.
  $$
   By~\cite[theorem 10]{paper}, $C_g$ is of simple type with $b_1=0$ and
   $b^+>1$. Alternatively, we can use~\cite{KM}, since 
   it contains a torus of self-intersection $0$ 
   intersecting an embedded $(-2)$-sphere transversely in one point
   (see proposition~\ref{prop:Cg}). 
\end{itemize}
\end{defn}

\begin{prop}
\label{prop:Bg}
  Consider $(B_g, \S_g)$. Let $K_{B_g}$ be the canonical class of 
  $B_g$, and $w=\PD [T_1] \in H^2(B_g;\ZZ)$.
  Then $\DD^w_{B_g}(e^{\bullet})={1 \over 2^{2g-2}} e^{Q/2} e^{K_{B_g}} +
  {1 \over 2^{2g-2}} e^{Q/2} e^{-K_{B_g}} + \cdots$, where the dots correspond
  to basic classes $\k$ for $B_g$ with $|\k \cdot \S_g| <2g-2$.
\end{prop}

\begin{pf}
  Write $B_g$ as $S_g \# g \overline{\CP}^2$. Let $F$ be the fibre of the
  natural elliptic fibration (i.e. $F=
  T_1$). Let $E_1, \ldots, E_g$ be the exceptional divisors. Then the
  basic classes are $k F \pm E_1 \pm E_2 \cdots \pm E_g$, with $-(g-2) \leq
  k \leq (g-2)$ and $k \equiv g-2 \pmod 2$ (see~\cite{KM}~\cite{FS}). 
  So the only basic class $\k$
  with $\k \cdot \S_g =2g-2$ is the canonical class $K_{B_g}=(g-2)F+ E_1
  +E_2 + \cdots +E_g$.  Therefore~\cite{KM}~\cite{FS}
$$   
   \DD_{B_g}(e^{\bullet})=e^{Q/2} (\sinh F)^{g-2} \sinh E_1 \cdots \sinh E_g
   ={1 \over 2^{2g-2}} e^{Q/2} e^{K_{B_g}} +
   {1 \over 2^{2g-2}} e^{Q/2} e^{-K_{B_g}} + \cdots
$$
  Now we note from proposition~\ref{prop:KM} that for $w=T_1$,
  $\DD^w_{B_g}=\DD_{B_g}$, since $w \cdot \k=0$ for all basic classes $\k$, 
  hence the result.
\end{pf}

\begin{prop}
\label{prop:Cg}
  $C_g$ is of simple type. For $w=\PD [\hat \S_2] \in H^2(C_g; \ZZ)$,
$$
  \DD^w_{C_g}(e^{\a})= 
  -2^{3g-5} e^{Q(\a)/2} e^{K_{C_g} \cdot \a} +
   (-1)^g 2^{3g-5} e^{Q(\a)/2} e^{K_{C_g} \cdot \a},
$$
  for a unique basic class $K=K_{C_g} \in H^2(C_g;\ZZ)$, such that
  $K\cdot \hat \S_2=2$, $K \cdot \S_g=2g-2$.
\end{prop}

\begin{pf}
By propositon~\ref{prop:Bg} and proposition~\ref{prop:DwS},
$$   
   D^{(T_2,\hat \S_2)}_{\hat B_2} (e^{\a})=
   {1 \over 4} e^{Q(\a)/2} e^{K_{\hat B_2} \cdot \a} +
   {1 \over 4} e^{Q(\a)/2} e^{-K_{\hat B_2} \cdot \a} + \cdots
$$
where the dots correspond to basic classes $\k$ with $\k \cdot \hat\S_2=0$.
Now we express
$C_g=\underbrace{\hat B_2 \#_{\hat \S_2} \cdots \#_{\hat \S_2} \hat B_2}_{g}$ 
and use~\cite[theorem 6]{paper}, 
$$
  D^{(\S_g,\hat \S_2)}_{C_g}(e^{\a})= 
  (-32)^{g-1} \left({1 \over 4}  \right)^g e^{Q(\a)/2} e^{K_{C_g} \cdot \a} +
  32^{g-1} \left({1 \over 4} \right)^g e^{Q(\a)/2} e^{K_{C_g} \cdot \a} =
$$
$$
  =(-1)^{g-1} 2^{3g-5} e^{Q(\a)/2} e^{K_{C_g} \cdot \a} +
   2^{3g-5} e^{Q(\a)/2} e^{K_{C_g} \cdot \a}, 
$$
where $K=K_{C_g} \in H^2(C_g;\ZZ)$ is defined as the unique cohomology
class such that
\begin{itemize}
\item $K \cdot \hat \S_2=2$
\item $K \cdot \S_g= 2g-2$
\item Writing $C_g= C_{g-1}^o \cup_Y \hat B^o_2$, one has 
      $K \cdot (\a_1 + \a_2)= K_{C_{g-1}} \cdot \a_1 +
      K_{\hat B_2} \cdot \a_2$, where  $\a_1 \in H_2(C_{g-1}^o)$
      and $\a_2 \in H_2(\hat B_2^o)$. 
\item Write $\hat B_2^o= (\tilde S_1- (N_{T_1} \cup N_{T_2})) 
  \cup_{Y-N_{T_1}\cap \bd N_{T_2}} 
  (\tilde S_1- (N_{T_1} \cup N_{T_2}))$. For any $\g \subset \S \subset
  \Y =\bd N_{T_1}$ disjoint from $N_{T_2}$, 
  we consider the vanishing discs for $\tilde S_1$ with elliptic 
  fibration with fibre $T_1$ (see~\cite[page 167]{Friedman}). These are
  embedded $(-1)$-discs. So, in the terminology of~\cite{MSz}, the preferred
  identification for $C_g= C_{g-1}^o \cup_Y \hat B^o_2$ is admissible. By
  remark~\ref{rem:adm},
  in the first place $C_g$ is of simple type. In the second place, there
  is a subspace $V \subset H_2(X)$ such that $H_2(X)=\cH \oplus V$ and
  $K \cdot \a=0$, for all $\a \in V$. 
\end{itemize}

For $w=\PD[\hat \S_2]$, using proposition~\ref{prop:DwS} again,
$$
  \DD^{\hat \S_2}_{C_g}(e^{\a})= 
  -2^{3g-5} e^{Q(\a)/2} e^{K_{C_g} \cdot \a} +
   (-1)^g 2^{3g-5} e^{Q(\a)/2} e^{K_{C_g} \cdot \a} 
$$
\end{pf}

Finally, we can find $M_{11}(t)$ and $M_{22}(t)$. Let $(X_1, \S_1)
=(X_2, \S_2)=(B_g,\S_g)$, $w_1=w_2=\PD[T_2]$, 
$D_1=D_2=T_2$ and $X=C_g$, $w=\PD[\hat \S_2]$, $D=\hat \S_2$, $\S=\S_g$, so
formula~\eqref{eqn:f5} and proposition~\ref{prop:Bg} yield
$$
  \DD^w_X(e^{tD+s\S})= e^{ts} (M_{11}( t ) {1 \over 2^{2g-2}}
    {1 \over 2^{2g-2}}
     e^{(2g-2)s} +  M_{22}( t ) {1 \over 2^{2g-2}}{1 \over 2^{2g-2}}
     e^{-(2g-2)s}) 
$$
and use proposition~\ref{prop:Cg} to get
$$
  \DD^w_X(e^{tD+s\S})= e^{ts}( -2^{3g-5} e^{2t+(2g-2)s} + 
  (-1)^g 2^{3g-5} e^{-2t-(2g-2)s}),
$$
 from where $M_{11}(t)=- 2^{7g-9} e^{2t}$, $M_{22}(t)=(-1)^g 2^{7g-9}
 e^{-2t}$.
 This finishes the proof of theorem~\ref{thm:main}.

\section{Proof of Theorem~\protect\ref{thm:finite}}

  Let $S= \S \x \CP^1$, $w= \PD [\CP^1] \in H^2(S;\ZZ)$. Suppose that
  $S$ is of finite type of order $n \geq 1$ with respect to $w$ and
  $w+\S$ for the metrics defined by $\S$. Then
  $$
   0=\Dws_S((x^2-4)^n z_{\a}e^{t\CP^1})=<\p^w (A,(x^2-4)^n e^{t\De}) , e_{\a}>
  $$
  for all $\a$. From lemma~\ref{lem:ea}, $ \p^w (A,(x^2-4)^n e^{t\De}) =0$.
  Then consider $(X,\S)$ permissible, $w \in H^2(X;\ZZ)$ with 
  $w \cdot \S \equiv 1 \pmod 2$, $X=\X \cup_Y A$.
  For any $D \in H_2(X)$ with $D \cdot \S=1$, we can write $D=\D +\De$,
  so 
  $$
   \Dws_X((x^2-4)^n e^{tD})=<\p^w (\X,e^{t\D}), \p^w (A,(x^2-4)^n e^{t\De})> 
   =0.
  $$ 
  We conclude that $\Dws_X((x^2-4)^n e^{tD})=0$ for all $D \in H_2(X)$,
  i.e. $X$ is of finite type 
  of order at most $n$ with respect to $w$.

\section{Conjecture}
\label{sec:conj}

Following our results in~\cite{paper} for the case of genus $g=2$, we
propose the following conjecture.

For any $(X,\S)$ permissible, define $\tilde X=X \#_{\S} B_g$ (we need to fix
one identification arbitrarily). Then $b^+(\tilde X)>1$ and also
$b_1(\tilde X)=0$ whenever $b_1(X)=0$. For any cycle $\D \subset
\X$ with $\bd \D=\g \subset Y$ an embedded curve
(when $\X$ has a cylindrical end,
$\D \cap (Y \x [0,\infty)) = \g\x [0,\infty)$ ), we choose cappings
$D=\D+\D_{B_g}$, $\bd \D_{B_g}=-\g$. The cappings have to satisfy
the following condition. In $B_g=S_g \# g
\overline{\CP}^2$ (see definition~\ref{def:Bg,Cg}), 
we fix an embedded surface $S$ representing $\s_g$,
intersecting $\S_g$ transversely in $g$ points. Put $S^o=S \cap \B$, then
we impose that $\D_{B_g} \cdot S^o=0$ (the pairing makes sense as long
as $\g$ and $\bd D^o$ are disjoint)

\begin{conj}
  Let $(X_1,\S_1)$ and $(X_2,\S_2)$ be permissible with $X_i$ having
  $b_1=0$ (we do not suppose that they are of simple
  type). Consider $\tilde X_i$. Then $\tilde X_i$ are of simple
  type. Put\/
  $\DD^w_{\tilde X_1}(e^{\a})=e^{Q(\a)/2}\sum \tilde a_{j,w} e^{\tilde 
  K_j\cdot \a}$ and\/
  $\DD^w_{\tilde X_2}(e^{\a})=e^{Q(\a)/2}\sum \tilde b_{k,w} e^{\tilde 
  L_k\cdot \a}$.
  Choose any identification $\p$ and let $X=X(\p)=X_1 \#_{\S} X_2$ 
  be the connected sum along $\S$, with the induded homology orientation.
  Then $X$ is of simple type.
  Choose $w \in H^2(X;\ZZ)$, $w_i \in H^2(\tilde X_i ;\ZZ)$ in a compatible
  way such that $w_i|_{B^o_g} = k \PD [T_1]|_{B^o_g}$.
  For every $D \in H_2(X)$, we write $D=\D_1 + \D_2$, $\D_i \subset \X_i$,
  with $\bd \D_1=
  - \bd \D_2=\g$ disjoint with $\bd S^o$
  and consider cappings $D_i \in H_2(\tilde X_i)$ of $\D_i$ as above, 
  in such a way that the map $D \mapsto (D_1,D_2)$ is linear. Then
\begin{eqnarray*}
   \DD^w_X(e^{tD}) &=& \f \,e^{Q(tD)/2}(\hspace{-5mm} \sum_{\tilde 
    K_j \cdot \S=\tilde L_k \cdot
    \S = 2g-2}\hspace{-5mm}  -2^{-3g+5} \tilde a_{j,w}\tilde b_{k,w} \, 
    e^{(\tilde K_j \cdot D_1 + 
    \tilde L_k \cdot D_2)t} + \\ 
    & & + \hspace{-5mm}\sum_{\tilde K_j \cdot \S=\tilde L_k \cdot \S= -(2g-2)}
    \hspace{-5mm} (-1)^g \,2^{-3g+5} \tilde a_{j,w}\tilde 
    b_{k,w} \, e^{(\tilde K_j \cdot
    D_1 +  \tilde L_k \cdot D_2)t}),
\end{eqnarray*}
for $\f=(-1)^{(g-1)(w^2-w_1^2-w_2^2)/2}$.
\end{conj}

\end{document}